\documentclass{article}
\pdfoutput=1

\usepackage{adjustbox}
\usepackage{subcaption}
\usepackage{hyperref}
\usepackage{algorithm}
\usepackage[noend]{algpseudocode}
\usepackage{amsmath}
\usepackage{amssymb}
\usepackage{amsthm}
\usepackage{tikz-cd}
\usepackage[capitalize]{cleveref}

\newcounter{Halgorithmic}
\AtBeginEnvironment{algorithmic}{\stepcounter{Halgorithmic}}
\ExpandArgs{c}\renewcommand{theHALG@line}{\arabic{Halgorithmic}.\arabic{ALG@line}}

\numberwithin{algorithm}{section}
\numberwithin{equation}{section}
\numberwithin{figure}{section}

\DeclareMathOperator{\divi}{div}
\DeclareMathOperator{\md}{~{\mathrm{mod}}~}
\DeclareMathOperator{\Vect}{Vect}
\DeclareMathOperator{\Set}{Set}
\DeclareMathOperator{\id}{id}

\newcommand{\tensorsum}[0]{\bigotimes\limits_{l = 1}^d\bigoplus\limits_{s_l \in [p_l]}}
\newcommand{\tensor}[0]{\bigotimes\limits_{l = 1}^d}
\newcommand{\sumtensor}[0]{\bigoplus\limits_{s \in p}\bigotimes\limits_{l = 1}^d}
\newcommand{\dsum}[0]{\bigoplus\limits_{s \in [p]}}
\newcommand{\prodsum}[0]{\prod\limits_{l = 1}^d\coprod\limits_{s_l \in [p_l]}}
\newcommand{\sumprod}[0]{\coprod\limits_{s \in p}\prod\limits_{l = 1}^d}

\pgfkeys{/pgf/number format/.cd, fixed, fixed zerofill, precision=0}

\newtheorem{theorem}{Theorem}

\theoremstyle{definition}
\newtheorem{defi}{Definition}

\title{Category Theory for Supercomputing: The Tensor Product of Linear BSP Algorithms}
\author{Thomas Koopman\thanks{Software Science,
                              Radboud University 
                              (\texttt{thomas.koopman@ru.nl})} 
        \and
        Rob H. Bisseling\thanks{Mathematical Institute,
                                Utrecht University,
                                (\texttt{R.H.Bisseling@uu.nl})}
        \and
        Sven-Bodo Scholz\thanks{Software Science,
                                Radboud University 
                                (\texttt{svenbodo.scholz@ru.nl})}}

\begin{document}
\maketitle

\begin{abstract}
We show that a particular class of parallel algorithm for linear
functions can be straightforwardly generalized to a parallel algorithm
of their tensor product. The central idea is to take a model of parallel
algorithms---Bulk Synchronous Parallel (\emph{BSP})---that decomposes parallel 
algorithms into so-called supersteps that are one of two types: a
\emph{computation superstep} that only does local computations, or a
\emph{communication superstep} that only communicates between processors.
We connect each type of supersteps to linear 
algebra with functors. Each superstep in isolation is simple enough
to compute their tensor product in Vect with well-known techniques of linear
algebra. We then individually translate the tensor product of supersteps back 
to the language of BSP algorithms. The functoriality of the tensor product 
allows us to compose the supersteps back into a BSP algorithm for the tensor
product of the original function. We state the recipe for creating these
new algorithms with only a minimal amount of algebra, so that it can be applied
without understanding the category-theoretic details.
We have previously used this to derive an efficient algorithm for the
higher-dimensional Discrete Fourier Transform, which we use as
an example throughout this paper. We also derive a parallel algorithm for
the Discrete Cosine Transform to demonstrate the generality of
our approach.
\end{abstract}

\section{Introduction}

Category theory is usually applied to other areas of abstract algebra,
or theoretical computer science. In this paper however, we look at a more
applied area: high-performance computing. It turns out that we can apply
central ideas behind category theory (compositionality and connecting different
areas of mathematics)
to the very practical problem of implementing a fast algorithm
for computing the Discrete
Fourier Transform (DFT) on a supercomputer. This algorithm---called
the Fast Fourier Transform (FFT)---implements the DFT in $O(n \log n)$ time
for $n$ data points and is considered one of the top 10 algorithms of the 20th
century~\cite{Dongarra00}.

Today, the `super' in the word `supercomputer' does not refer to a single processor that can compute
very quickly, but rather to the number of processors. In fact, a
supercomputer is a cluster of many connected processors, each with their own 
memory. Accessing data on a different processor requires explicit data movement 
and minimising this while maintaining load balance
is the essence of making programs on supercomputers faster.
We can reason about this by separating computation from data movement with a
computational model called the 
\emph{Bulk Synchronous Parallel} model~\cite{Valiant90}, or
\emph{BSP} for short.

More specifically, the contribution of this paper is \cref{thm:tensor-of-bsp},
a recipe for taking the tensor product of certain BSP algorithms, formulated
without category theory. The key idea behind the proof, however,
is to use category theory as a tool to transfer
knowledge from one area of mathematics to another, in this case parallel
algorithms and linear algebra. We can connect a certain
class of BSP algorithms to constructs in linear algebra. This uncovers structure
that is well-understood in the category Vect, and we use this theory to derive the core
of our algorithm, namely distributing the tensor product over the direct sum.
The compositional nature of BSP algorithms allows us to study simple steps in
isolation and combine the result.
The algebraic result can then be translated back to the category Set, giving a BSP algorithm of the same structure for
higher-dimensional arrays. The application of \cref{thm:tensor-of-bsp} to
the one-dimensional algorithm of~\cite{Inda01} has previously been
published~\cite{Koopman23} and benchmarked, showing improvements over the 
state-of-the-art. Though that paper provides an elementary proof, the algorithm 
was derived using the theory from the present paper.

\subsection{The Discrete Fourier Transform}

The discrete Fourier Transform $F_n$ of length $n$ is a linear function 
$\mathbb{C}^n \to \mathbb{C}^n$. Writing $\omega_n$ for the $n$th root of unity 
$e^{-2\pi i / n}$, it is given explicitly by
\[
    y_k = \sum_{j = 0}^{n - 1} x_j \omega_n^{jk}.
\]
There is also a multidimensional variant of the DFT that operates
on $n_1 \times \cdots \times n_d$ arrays:
\begin{equation}\label{eq:higher-rank-DFT}
    y_{k_1,  \cdots , k_d} = \sum_{j_1 = 0}^{n_1 - 1}  \cdots  \sum_{j_d = 0}^{n_d - 1} x_{j_1,  \cdots , j_d} \omega_{n_1}^{j_1k_1} \cdot  \cdots  \cdot \omega_{n_d}^{j_dk_d}.
\end{equation}
We will call an $n_1 \times \cdots \times n_d$ array, an array of rank $d$ 
rather than dimension $d$ from now on in the tradition of the programming language
APL~\cite{Iverson62}.
For notational convenience, we will define $[n] := \{0, \cdots, n - 1\}$
and abbreviate $(k_1, \cdots, k_d) \in [n_1] \times \cdots \times [n_d]$
as $k$.

The original motivation behind this work is to find an efficient parallel
algorithm for the multidimensional DFT. To this end, it is important to
understand what \cref{eq:higher-rank-DFT} means algebraically.
The central concept behind generalizing the DFT from rank-1 to higher ranks is multilinearity:
functions $f: V_1 \times \cdots \times V_d \to W$ that are linear in each
argument. We may ask ourselves: what is the vector space that best captures 
multilinearity from $V_1 \times \cdots \times V_d$? The answer is the tensor
product $V_1 \otimes \cdots \otimes V_d$, which together with a multilinear
map $\otimes: V_1 \times \cdots \times V_d \to
V_1 \otimes \cdots \otimes V_d$ can be characterised by the universal property
that for any multilinear map $f: V_1 \times \cdots \times V_d \to W$,
there exists a unique linear map 
$\tilde{f}: V_1 \otimes \cdots \otimes V_d \to W$ such that 
$\tilde{f} \circ \otimes = f$.

\begin{center}
\begin{tikzcd}[column sep = huge, row sep = huge]
    V_1 \times \cdots \times V_d \arrow[d, "\otimes"'] \arrow[r, "f"] & W \\
    V_1 \otimes \cdots \otimes V_d \arrow[ur, dotted, "\tilde{f}"'] &
\end{tikzcd}
\end{center}

We can use this to turn a product of functions 
$(f_1, \cdots, f_d): V_1 \times \cdots \times V_d \to 
W_1 \times \cdots \times W_d$ into a tensor product 
$f_1 \otimes \cdots \otimes f_d: V_1 \otimes \cdots \otimes V_d \to 
W_1 \otimes \cdots \otimes W_d$ by traversing the following commutative diagram
from bottom left to top right to bottom right. Working out the tensor product
of rank-1 DFTs yields the higher-ranked DFT of \cref{eq:higher-rank-DFT}.

\begin{center}
\begin{tikzcd}[column sep = huge, row sep = huge]
    V_1 \times \cdots \times V_d \arrow[d, "\otimes"'] \arrow[r, "{(f_1, \cdots, f_d)}"] & W_1 \times \cdots \times W_d \arrow[d, "\otimes"] \\
    V_1 \otimes \cdots \otimes V_d \arrow[ur, dotted] & W_1 \otimes \cdots \otimes W_d
\end{tikzcd}
\end{center}

\subsection{Distributed Computing}

\begin{figure}
    \begin{subfigure}{0.55\textwidth}
        \begin{tikzpicture}[scale=0.70]
            \foreach \j in {0, ..., 2} {
                \draw[fill=red] (0 + 3 * \j, 0) rectangle (0 + 3 * \j + 1, 1);
            }
            \foreach \j in {0, ..., 2} {
                \draw[fill=yellow] (1 + 3 * \j, 0) rectangle (1 + 3 * \j + 1, 1);
            }
            \foreach \j in {0, ..., 2} {
                \draw[fill=green] (2 + 3 * \j, 0) rectangle (2 + 3 * \j + 1, 1);
            }
            \foreach \j in {0, ..., 8} {
                \node at (\j + 0.5, 0.5) {$f(\j)$};
            }
        \end{tikzpicture}
        \caption{Global view}
    \end{subfigure}
    \hfill
    \begin{subfigure}{0.4\textwidth}
        \begin{tikzpicture}[scale=0.75]
            \foreach \j [evaluate=\j as \jeval using 3 * \j + 2] in {0, ..., 2} {
                \draw[fill=green] (\j, -1.5) rectangle 
                        node {$f(\pgfmathprintnumber{\jeval})$} (\j + 1, -1.5 + 1);
            }
            \draw node at (4.5, -1) {$\color{green}{P(2)}$};
            \foreach \j [evaluate=\j as \jeval using 3 * \j + 1] in {0, ..., 2} {
                \draw[fill=yellow] (\j, 0) rectangle 
                        node {$f(\pgfmathprintnumber{\jeval})$} (\j + 1, 0 + 1);
            }
            \draw node at (4.5, 0.5) {$\color{orange}{P(1)}$};
            \foreach \j [evaluate=\j as \jeval using 3 * \j + 0] in {0, ..., 2} {
                \draw[fill=red] (\j, 1.5) rectangle 
                        node {$f(\pgfmathprintnumber{\jeval})$} (\j + 1, 1.5 + 1);
            }
            \draw node at (4.5, 2) {$\color{red}{P(0)}$};
        \end{tikzpicture}
        \caption{Local view}
    \end{subfigure}
    \caption{Graphical representation of cyclically distributed initialisation 
             over three processors, indicated in red, yellow, and green.}
    \label{fig:cyclic}
\end{figure}
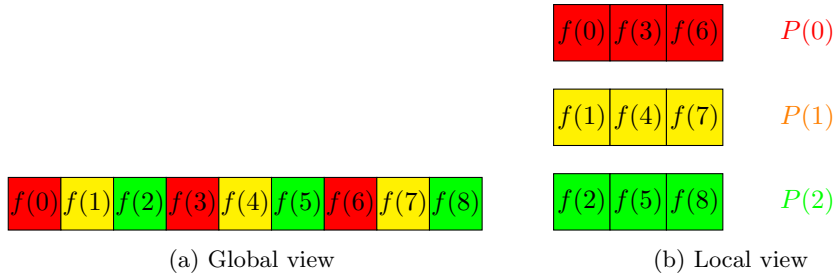

In programming, there is usually little difference between data structures
and the mathematical objects they represent. In the programming language
Fortran, we would write
\texttt{X(k)} for $x_k$ or \texttt{X(k1, k2)} for $x_{k_1, k_2}$. The 
programming language can then easily compute the address in memory from such an 
expression. This approach fails when the data structures are so large that they
do not fit in the memory of a single machine. Programming a supercomputer
is also called distributed-memory parallel computing because we have to distribute the data
structure over different processors with their own memory. If we have
$p$ processors, we will identify processors with an index $0 \leq s < p$.
Each processor $P(s)$ has its local data structure $X^{(s)}$. These local data structures then
implicitly define the global data structure $X$ through a correspondence
called the \emph{distribution}. \Cref{alg:cyclic} gives an
example, called the \emph{cyclic distribution}. In this paper, the \textbf{to} in
algorithms is always exclusive. We can view $X$ as one data structure, called
the \emph{global view}, or as a collection of local data structures $(X^{(s)})_{s}$,
called the \emph{local view}, as illustrated in \cref{fig:cyclic}. The
distribution can be interpreted as a bijection between the two views.

\begin{figure}
    \begin{subfigure}{0.45\textwidth}
        \begin{algorithmic}
            \For{$j = 0$ \textbf{to} $n$}
                \State $X[j] = f(j)$
            \EndFor
        \end{algorithmic}
        \caption{Sequential initialisation}
    \end{subfigure}
    \quad
    \begin{subfigure}{0.45\textwidth}
        \begin{algorithmic}
            \For{$k = 0$ \textbf{to} $n / p$}
                \State $X^{(s)}[k] = f(s + kp)$
            \EndFor
        \end{algorithmic}
        \caption{Distributed initialisation on $P(s)$}
    \end{subfigure}
    \caption{Initialisation of $x_j = f(j)$ sequentially, and in parallel.}
    \label{alg:cyclic}
\end{figure}

\subsection{The Bulk Synchronous Parallel Model}

The Bulk Synchronous Parallel (BSP) model is often used, either implicitly
or explicitly, for parallel
computing. It comprises a parallel computer architecture, a class of parallel 
algorithms, and a function for charging costs to algorithms~\cite{Valiant90}. 
We employ the variant extensively used in~\cite{Bisseling20} and focus on the algorithmic part. 

\emph{A BSP algorithm} consists of a sequence of supersteps. We have
two kinds of supersteps. The first, a \emph{computation} superstep,
performs a sequence of operations on local data structures $X^{(s)}$.
In the second, a \emph{communication} superstep, each processor sends and 
receives a number of messages. At the end of a superstep, all processors 
synchronise as follows. Each processor checks whether it has finished all its 
tasks of that superstep. Processors wait until all others have 
finished. When this has happened, they all proceed to the next 
superstep. This form of synchronisation is called \emph{bulk synchronisation}, 
hence the name of the model.

\section{A BSP Algorithm for the Rank-1 DFT}

To get a feel for BSP algorithms, we will derive a BSP algorithm for
the four-step FFT framework, which can be
found in a slightly different formulation in~\cite{Loan92}. We use different
variable names and notation to make it easier to understand the parallel 
algorithm. To that end, suppose $x$ is an array of length $n$, that $p \mid n$ and 
$p \mid \frac{n}{p}$ (or equivalently that $p^2 \mid n$), and that we want to 
calculate $y := F_n(x)$. We write $v(a:b:c)$ to mean the strided subarray of 
$v$ which starts at index $a$ and has stride $b$. The last variable $c$ is the 
length of $v$. If we write $v(a:b)$ we mean the subarray that starts at $a$ and 
ends at $b$ (exclusive). This lets us write the four-step framework as the sequential 
\cref{chp1:four step sequential}.

\begin{algorithm}[ht]
\caption{Sequential four-step FFT framework}\label{chp1:four step sequential}
\begin{algorithmic}[1]
    \Require $x: $ array of length $n$, number $p$ such that $p^2 \mid n$.
    \Ensure $y: $ array of length $n$, such that $y = F_n(x)$.
    \For{$s \in [p]$} \Comment{Step 1}
		\State $x^{(s)} = x(s: p: n)$;
		\State $x^{(s)} = F_{n/p}(x^{(s)})$;
	\EndFor
	
    \For{$s \in [p]$} \Comment{Step 2}
        \For{$k \in [n/p]$}
			\State $x^{(s)}_k = \omega_n^{ks} x^{(s)}_k$;
		\EndFor
	\EndFor
	
    \For{$s \in [p]$} \Comment{Step 3}
        \For{$k \in [n/p]$}
			\State $y(k: \frac{n}{p}: n)_s = x^{(s)}_k$;
		\EndFor
	\EndFor
		
    \For{$k \in [n/p]$} \Comment{Step 4}
		\State $y(k: \frac{n}{p}: n) = F_{p}(y(k : \frac{n}{p} : n))$;
	\EndFor
\end{algorithmic}
\end{algorithm}

We use the parallelisation strategy of \cite{Inda01}. 
The notation we chose in \cref{chp1:four step sequential} already
suggests using the cyclic distribution. This gives us the first two steps of 
parallel \cref{chp1:four step parallel}. A BSP algorithm in Single Program,
Multiple Data (SPMD) style is written 
from the perspective of a processor, so the loop over $s$ is removed
and $s$ becomes the processor index. 

\begin{algorithm}[ht]
\caption{Parallel four-step FFT framework for $P(s)$}\label{chp1:four step parallel}
\begin{algorithmic}[1]
    \Require $x: $ array of length $n$, distr$(x) = $ cyclic over $p$ processors such that $p^2 \mid n$.
    \Ensure $y: $ array of length $n$, distr$(y) = $ cyclic, such that $y = F_n(x)$.
	\State $x^{(s)} = F_{n/p}(x^{(s)})$; \Comment{Step 1}
    \For{$k \in [n/p]$} \Comment{Step 2}
			\State $x^{(s)}_k = \omega_n^{ks} x^{(s)}_k$;
	\EndFor
	
    \For{$k \in [p]$}\Comment{Step 3}
		\State Put $x^{(s)}(k : p : \frac{n}{p})$ in $P(k)$ as $y^{(k)}(s \frac{n}{p^2}: (s + 1)\frac{n}{p^2})$;
	\EndFor
	
    \For{$t \in [n/p^2]$}\Comment{Step 4}
		\State $y^{(s)}(t : \frac{n}{p^2} : \frac{n}{p}) = F_p(y^{(s)}(t : \frac{n}{p^2} : \frac{n}{p}))$;
	\EndFor
\end{algorithmic}
\end{algorithm}

In the BSP framework we view Step 1 and Step 2 together as one computation
superstep, Step 3 as a communication superstep, and Step 4 as a computation
superstep as well. Steps 3 and 4 are derived as follows.

To parallelize Step 3, we decompose $k$ as
$k = k' + cp, \quad 0 \leq k' < p, \quad 0 \leq c < n/p^2$ so we
can consider the $k$th element of $x^{(s)}$ to be the $c$th element of the strided
subarray $x^{(s)}(k' : p : n / p)$. The following index chase shows that
this strided subarray becomes the $s$th contiguous subarray on $P(k')$.
\begin{align*}
    y(k : n / p : n)_s &= y_{k + sn/p} \\
                       &= y_{k' + (sn/p^2 + c)p} \\
                       &= y^{(k')}_{sn/p^2 + c} \\
                       &= y^{(k')}(sn/p^2 : (s + 1)n/p^2)_c.
\end{align*}
Finally, we drop the prime from $k'$ to obtain Step 3 of the parallel
algorithm.

To parallelize Step 4, we first decompose $k$ as
$k = tp + s, \quad 0 \leq s < p, \quad 0 \leq t < n/p^2$. We then use our 
assumption that $p^2 \mid n$ to rewrite the $c$th index of $y(k : n/p : n)$ as
$k + c\frac{n}{p} = s + (t + cn/p^2)p$. We conclude that the $c$th element
of $y(k : n/p : n)$ is the $(t + cn/p^2)$th element of $y^{(s)}$,
which means it is the $c$th element of $y^{(s)}(t : n/p^2 : n / p)$.

As a result, we have obtained a complete BSP algorithm for the fast computation of a rank-1 DFT.

\section{Linear BSP Algorithms}

To apply theory from linear algebra, we must first translate the structure
of a BSP algorithm to Vect. We restrict the class of BSP algorithms
to those that compute linear functions using supersteps that respect
linearity in a way that we define more precisely in \cref{defi:linear-bsp}.
The BSP model does not explicitly mention distributions because they
do not incur computation or data movement, but they are algebraically
important. We discuss distributions, computation supersteps, and communication
supersteps.

\subsection{Arrays as Vectors}

A computer program for a linear function $f$, is not bound by linear structure.
So writing $U$ for the forgetful functor, we compute $Uf$.

The forgetful functor $U: \Vect \to \Set$ is accompanied by a free functor
$V: \Set \to \Vect$ that adds structure. This functor is used whenever we
work with arrays: an array $X$ corresponds to the vector
$\sum_{j \in [n]} X[j] b_j$, where $(b_j)_{j \in [n]}$ is a basis of
$\mathbb{C}^n$. The free functor applied to a set $B$ returns a vector space
$VB$ consisting of all formal linear combinations of elements in $B$ (so
$B$ is a basis of $VB$). This means that we can interpret the set of
arrays over an index set $I$ as $UVI$. We can use this view to explain why we 
represent the tensor product of DFTs as an operation on 
$n_1 \times \cdots \times n_d$ arrays in \cref{eq:higher-rank-DFT}.
We can form a basis of the tensor
product $VB_1 \otimes \cdots \otimes VB_d$ by taking the product of bases
$B_1 \times \cdots \times B_d$, so $VB_1 \otimes \cdots \otimes VB_d \cong
V(B_1 \times \cdots \times B_d)$. Given rank-$1$ index sets of the form 
$I = [n]$, this gives a $d$-dimensional index set of 
the form $[n_1] \times \cdots \times [n_d]$, which indexes a rank-$d$ array.

\subsection{Distributions}

In the BSP setting we have a collection of index sets 
$(I^{(s)})_{s \in [p]}$ for the processors, rather than one global 
index set $I$. A distribution is then a collection of functions 
$\phi^s: I^{(s)} \to I$ that links these through $X^{(s)}[k] = X[\phi^s(k)]$.
In other words, we have a disjoint union 
$\coprod_{s \in [p]}I^{(s)}$ and an induced function
$\phi = [\phi^s \mid s \in [p]]: \coprod_{s \in [p]} I^{(s)} \to I$
that is bijective. The free functor turns this into the bijection between the 
local and global view of \cref{fig:cyclic}.

Viewing distributions as bijections between index sets allows for an algebraic 
interpretation of parallelism: the disjoint union
is the coproduct in Set, and as $V \dashv U$, the free functor preserves this
colimit and maps it to the coproduct in Vect, the direct sum. The direct sum
is also the product in Vect, which is preserved by the forgetful functor.
This gives the commutative diagram of \cref{diag:distr}, where we use the square
braces for the unique morphism out of coproducts in both Set and Vect.

\begin{figure}[ht]
    \centering
\begin{tikzcd}[row sep = huge]
    UV\coprod\limits_{s \in [p]} I^{(s)} \arrow[r, "\cong"] \arrow[d, "{UV[\phi^s \mid s \in [p]]}"] & 
        U\dsum VI^{(s)} \arrow[r, "\cong"] \arrow[d, "{U[V\phi^s \mid s \in [p]]}"] &
            \prod\limits_{s \in [p]} UVI^{(s)} \\
    UVI \arrow[r, "="] & UVI \arrow[r, "="] & UVI
\end{tikzcd}
    \caption{Distributions correspond to the direct sum in Vect. Local view at
             the top, global view at the bottom.}
    \label{diag:distr}
\end{figure}

The left side of \cref{diag:distr}, shows that a distributed array is really
just an array with index set $\coprod_{s \in [p]} I^{(s)}$, and that the
distribution is an isomorphism between this array and $X$. In Fortran, this is how
we think about distributed computing when using coarrays:
we can declare a distributed array of rank two as \texttt{X(:,:)[*]}, and
then refer to $X^{(s)}$ by \texttt{X(:,:)[s]}.
The top right side, $\prod_{s \in [p]}UVI^{(s)}$, views a distributed array as
a collection of local data structures. This is how we usually think about
distributed computing in the BSP model. The top middle part exposes the algebraic
meaning of parallelism: the direct sum.

\subsection{Computation Superstep}

A computation superstep is a collection of functions between these data
structures, so a product. If these functions are linear, the product is
preserved under the forgetful functor, so they correspond to a direct sum
as illustrated in \cref{diag:bsp_comp}. Note that this holds for
\cref{chp1:four step parallel} as $F_{n/p}$, $F_p$ and
$x_k^{(s)} \mapsto \omega_n^{ks} x_k^{(s)}$ are all linear.

\begin{figure}[ht]
    \centering
    \begin{tikzcd}[row sep = huge]
        \prod\limits_{s \in [p]} UVI^{(s)} \arrow[r, "\cong"] \arrow[d, "\prod\limits_{s \in [p]} Uf^{(s)}"'] & 
            U\dsum VI^{(s)} \arrow[d, "U\dsum f^{(s)}"] \\
        \prod\limits_{s \in [p]} UVI^{(s)} \arrow[r, "\cong"]  & 
            U\dsum VI^{(s)} \\
    \end{tikzcd}
    \caption{A BSP computation superstep of linear functions}
    \label{diag:bsp_comp}
\end{figure}

\subsection{Communication Superstep}

The communication superstep of \cref{chp1:four step parallel} permutes
the local arrays, meaning it can be described as an isomorphism on the
index set $\coprod_{s \in [p]} I^{(s)}$. In our experience, many communication
supersteps in parallel algorithms can be described this way, so we restrict
ourselves to communication supersteps of this type.

As functions out of a coproduct are determined uniquely by their components,
we can write such an isomorphism as $r = [r_s \mid s \in [p]]$.
Concretely, we have $r_s(k) = (k \md p, sn/p^2 + k \divi p)$
in Step 3 of \cref{chp1:four step parallel} because we put the $k$th element of $P(s)$
in the $(sn/p^2 + k \divi p)$th position on $P(k \md p)$.

A permutation on basis elements is lifted to an isomorphism between vector
spaces by the free functor, so we consider isomorphisms as in
\cref{diag:bsp_comm}. 

\begin{figure}[ht]
    \centering
    \begin{tikzcd}[row sep = huge]
    UV\coprod\limits_{s \in [p]} I^{(s)} \arrow[d, "{UV[r_s \mid s \in [p]]}"]  \\
    UV\coprod\limits_{s \in [p]} I^{(s)}
\end{tikzcd}
    \caption{A BSP communication superstep that redistributes data 
             is a bijection on the index set under the free functor}
    \label{diag:bsp_comm}
\end{figure}

\subsection{Definition}

To conclude, the type of BSP algorithm that we consider factorises a linear
function $f$ in direct sums of functions and index permutations under the
free functor. We then forget the linear structure by using the forgetful
functor. More precisely:

\begin{defi}[Linear BSP algorithm]\label{defi:linear-bsp}
    Interpret arrays $X$ and $Y$ over index sets $I_x$, $I_y$ as $Ux$, $Uy$
    for vectors $x \in VI_x$, $y \in VI_y$. A \emph{linear BSP algorithm}
    has an array $X$ as input, an array $Y$ as output, and computes
    $y = f(x)$ for a linear function $f$ such that the following conditions
    hold:
    \begin{enumerate}
        \item all computation supersteps are of the form
              $Y^{(s)} = Uf^{(s)}(X^{(s)})$ for linear functions $f^{(s)}$;
        \item all communication supersteps are of the form $UVr$ where
              $r$ permutes the index set $\coprod_{s \in [p]}I^{(s)}$.
    \end{enumerate}
\end{defi}

\section{Tensor Product of Linear BSP Algorithms}

We can now state \cref{thm:tensor-of-bsp}, which shows that if we have linear
BSP algorithms for functions $f_1, \cdots, f_d$, we can straightforwardly derive
a BSP algorithm of the same structure for their tensor product
$f_1 \otimes \cdots \otimes f_d$.

\begin{theorem}[The tensor product of linear BSP algorithms]
\label{thm:tensor-of-bsp}
    Linear BSP algorithms $\mathcal{A}_1, \cdots \mathcal{A}_d$ for
    $f_1, \cdots, f_d$ with the same number and structure of supersteps can be combined into
    a linear BSP algorithm for $f_1 \otimes \cdots \otimes f_d$ by combining
    the supersteps and distributions as follows.

    \begin{itemize}
        \item If $\mathcal{A}_l$ uses $p_l$ processors indexed between
              $0$ and $p_l$, the higher-dimensional algorithm uses a
              $d$-dimensional processor grid
              $p := [p_1] \times \cdots \times [p_d]$.
        \item If $\mathcal{A}_l$ has distribution
              \[
                  [\phi^{(s_l)} \mid s_l \in [p_l]],
              \]
              the higher-dimensional algorithm has distribution
              \[
                  [\phi^{s_1} \times \cdots \times \phi^{s_d} 
                        \mid s \in p].
              \]
        \item Computation supersteps that compute $Ug_1, \cdots, Ug_d$ are combined into a
              computation superstep that computes
              $U(g_1 \otimes \cdots \otimes g_d)$.
        \item Communication supersteps that perform a permutation
              \[
                  Z^{(\psi_l(s_l, k_l))}[\rho_l(s_l, k_l)] = X^{(s_l)}[k_l]
              \]
              are combined into a communication superstep that performs the
              permutation
              \[
                  Z^{((\psi_1 \times \cdots \times \psi_d)(s, k))}
                   [(\rho_1 \times \cdots \times \rho_d)(s, k)] =
                   X^{(s)}[k].
              \]
    \end{itemize}
\end{theorem}
(Note that we can always add identity functions to make the number and structure
of supersteps match.)

The proof relies on two facts. First, the tensor \textbf{product} and direct 
\textbf{sum} distribute over each other in the same way the product and sum of 
numbers do:
\begin{equation}\label{eq:tensor-dsum}
    \left(\bigoplus_{s_1 \in [p_1]} VI^{(s_1)}\right) \otimes \cdots \otimes \left(\bigoplus_{s_d \in [p_d]} VI^{(s_d)}\right) \cong
    \bigoplus_{s \in p} (VI^{(s_1)} \otimes \cdots \otimes VI^{(s_d)}).
\end{equation}
This isomorphism is a natural transformation and allows us to pull the direct 
sum (encoding parallelism) out of the tensor product. 
Second, the tensor product is functorial, so we can study
the distribution, computation superstep, and communication superstep in 
isolation.

\subsection{Distributions}

\Cref{eq:tensor-dsum} is the intuition behind the tensor product turning a 
collection of distributed rank-$1$ arrays into a distribution of a rank-$d$ 
array, with also a rank-$d$ processor grid. We work this out in more detail in 
\cref{diag:distr:tensor}, where we also use that the tensor product is the 
image of the Cartesian product in Set under the free functor. We can tell that 
the $\phi_{s_l}$ are combined dimension-wise by comparing the right side of 
\cref{diag:distr:tensor} with the left side of \cref{diag:distr}.

\begin{figure}[ht]
    \centering
    \begin{adjustbox}{max width=\textwidth}
\begin{tikzcd}[row sep = huge]
    \tensorsum VI^{(s_l)} \arrow[r, "\cong"] \arrow[d, "{\tensor[ V \phi^{s_l}]}"] &
      \sumtensor VI^{(s_l)} \arrow[r, "\cong"] \arrow[d, "{[\tensor V \phi^{s_l}]}"] &
       \dsum V \prod\limits_{l = 1}^d I^{(s)}_l \arrow[d, "{[V\prod\limits_{l = 1}^d \phi^{s_l}]}"] \arrow[r, "\cong"] &
          V \coprod_{s = 0}^{p - 1} \prod\limits_{l = 1}^d I^{(s)}_l \arrow[d, "{V[\prod\limits_{l = 1}^d \phi^{s_l}]}"]\\
    \tensor V I_l \arrow[r, "="] & \tensor V I_l \arrow[r, "\cong"] & V\prod\limits_{l = 1}^d I_l \arrow[r, "="] & V\prod\limits_{l = 1}^d I_l
\end{tikzcd}
    \end{adjustbox}
    \caption{The tensor product of a distribution is again a distribution,
             where we take the Cartesian product of index sets and distribution
             functions}
    \label{diag:distr:tensor}
\end{figure}

\subsection{Computation Superstep}

\Cref{eq:tensor-dsum} directly gives the diagram of 
\cref{diag:bsp_comp:tensor}. We conclude that we can simply apply the tensor 
product locally.

\begin{figure}[ht]
    \centering
    \begin{tikzcd}[row sep = huge]
        \tensorsum VI_{d}^{(s_l)} \arrow[d, "\tensorsum f^{(s_l)}"'] \arrow[r, "\cong"] &
            \sumtensor VI_{d}^{(s_l)} \arrow[d, "\sumtensor f^{(s_l)}"] \\
        \tensorsum VI_{d}^{(s_l)} \arrow[r, "\cong"] &
            \sumtensor VI_{d}^{(s_l)} \\
    \end{tikzcd}
    \caption{The tensor product of a BSP computation superstep of linear 
             functions}
    \label{diag:bsp_comp:tensor}
\end{figure}

\subsection{Communication Superstep}

The right side of \Cref{diag:bsp_comm:tensor} shows that the tensor
product of a communication superstep is the product of the underlying
permutations composed with an isomorphism. 

\begin{figure}[ht]
    \centering
    \begin{tikzcd}[row sep = huge, column sep = huge]
        \tensor V \coprod\limits_{s_l \in [p_l]} I^{(s_l)}
            \arrow[d, "{\tensor V [r_{s_l} \mid s_l \in [p_l]]}"'] \arrow[r, "\cong"] &
            V \prodsum I^{(s_l)} \arrow[d, "{V \prod\limits_{l = 1}^d [r_{s_l} \mid s_l \in [p_l]]}"']  \arrow[r, "\cong"] &
                V \sumprod I^{(s_l)} \arrow[dl, "{V[\prod\limits_{l = 1}^d r_{(s_l)} \mid s \in p]}"] \\
        \tensor V \coprod\limits_{s_l \in [p_l]} I^{(s_l)} \arrow[r, "\cong"'] &
            V \prodsum I^{(s_l)} \arrow[r, "\cong"'] &
                V \sumprod I^{(s_l)} \\
    \end{tikzcd}
    \caption{The tensor product of a redistribution is done dimension-wise}
    \label{diag:bsp_comm:tensor}
\end{figure}

To show how this corresponds to \cref{thm:tensor-of-bsp}, we chase through the 
right side of the diagram. We write 
$r_{s_l}(k) = (\psi_l(s_l, k), \rho_l(s_l, k))$, which is interpreted as: 
Put $X^{(s_l)}[k]$ into $P(\psi_l(s_l, k))$ at local 
index $\rho_l(s_l, k)$. Starting at the top right corner, and ignoring
the free functor $V$, we get

\begin{align*}
    &  ((s_1, \cdots, s_d), (k_1, \cdots, k_d)) \mapsto \\
    &          (r_{s_1}(k_1), \cdots, (r_{s_d}(k_d)) = \\
    &        ((\psi_1(s_1, k_1), \rho_1(s_1, k_1)), \cdots, 
                (\psi_d(s_d, k_d), \rho_d(s_d, k_d)) \mapsto \\
    &        ((\psi_1(s_1, k_1), \cdots, \psi_d(s_d, k_d)),
                      (\rho_1(s_1, k_1), \cdots, \rho_d(s_d, k_d)).
\end{align*}

\section{Applications}

\subsection{Discrete Fourier Transform}

\Cref{thm:tensor-of-bsp} can be used to directly derive 
\cref{chp1:four step parallel:tensor} from \cref{chp1:four step parallel},
without the tedious proof of~\cite{Koopman23}.

\begin{algorithm}[ht]
\caption{Parallel four-step FFT framework for processor $P(s)=P(s_1, \cdots, s_d)$,
         rank-$d$}\label{chp1:four step parallel:tensor}
\begin{algorithmic}[1]
\Require $X: $ array of size $n_1 \times \cdots \times n_d$, distr$(X) = $ 
         rank-$d$ cyclic over $p_1 \times \cdots \times p_d$ processors such 
         that $p_l^2 | n_l$, for $l=1,\cdots,d$.
\Ensure $Y: $ array of size $n_1 \times \cdots \times n_d$, 
        distr$(Y) = $  rank-$d$ cyclic, such that 
        $Y = (F_{n_1} \otimes \cdots \otimes F_{n_d})(X)$.
	\State $X^{(s)} := (F_{n_1/p_1} \otimes \cdots \otimes F_{n_d/p_d})(X^{(s)})$; \Comment{Superstep 0}
	\For{$k \in [n_1/p_1] \times \cdots \times [n_d/p_d]$}
		\State $X^{(s)}[k] := (\prod_{l = 1}^d\omega_{n_l}^{k_ls_l}) X^{(s)}[k]$;
	\EndFor
	\For{$k \in [p_1] \times \cdots \times [p_d]$}\Comment{Superstep 1}
		\State Put $X^{(s)}(k : p : \frac{n}{p})$ in $P(k)$ as $Y^{(k)}[s\frac{n}{p^2}: (s + 1)\frac{n}{p^2}]$;
	\EndFor
	\For{$t \in [n_1/p_1^2] \times \cdots \times [n_d/p_d^2]$} \Comment{Superstep 2}
		\State $Y^{(s)}(t:\frac{n}{p^2}: \frac{n}{p}) := (F_{p_1} \otimes \cdots \otimes F_{p_d})\left(Y^{(s)}(t:\frac{n}{p^2}: \frac{n}{p})\right)$;
	\EndFor
\end{algorithmic}
\end{algorithm}

\subsection{Discrete Cosine Transform}

We can also apply \cref{thm:tensor-of-bsp} to the Discrete Cosine Transform
(DCT-II), see e.g.~\cite{Loan92,Press07}:
\begin{equation}\label{eq:dct-ii}
    y_k = \sum_{j = 0}^{n - 1} x_j \cos \frac{(2j + 1) k\pi}{2n}.
\end{equation}

The most efficient sequential rank-1 DCT-II algorithms pack the $n$ real input
points into $n/2$ complex data points, apply the DFT, and then extract the
result. Our approach is not applicable to these algorithms because the
extraction is linear with respect to $\mathbb{R}$, but not $\mathbb{C}$. This
mix of ground fields breaks the theory we have developed.
We discuss this further in \cref{subsec:future-tensor}.
Instead we use a less efficient algorithm that views the DCT-II as a complex
function, that we happen to apply to real input only.

We can implement the DCT-II by extending the signal, applying a DFT, and
then extracting the result~\cite{Makhoul80}. Given an input $x$ of length $n$, 
we define $w$ 
of length $2n$ as
\[
    w_j = \begin{cases}
        x_j & 0 \leq j < n \\
        x_{2n - 1 - j} & n \leq j < 2n.
    \end{cases}
\]
After transforming $z = DFT(w)$, the DCT-II $y$ can be extracted by
\[
    y_k = \frac{1}{2}\omega_{2n}^{k / 2} z_k, \quad 0 \leq k < n.
\]

\subsubsection{Linear BSP Algorithm}

To derive a linear BSP algorithm, we first decompose this algorithm into
functions. We write $r$ for the reversal function $x_j \mapsto y_{n - 1 - j}$,
$\langle f, g \rangle$ for the function $V \to W_1 \oplus W_2$ induced by
$f: V \to W_1$, $g: V \to W_2$, $\pi_1$ for the projection to the first
component, and $\mathbb{C}^n \oplus \mathbb{C}^n \cong \mathbb{C}^{2n}$ for
putting the bases after each other. This gives the following decomposition.

\begin{center}
    \begin{tikzcd}
        \mathbb{C}^n \arrow[r, "{\langle \id, \id \rangle}"] & 
        \mathbb{C}^n \oplus \mathbb{C}^n \arrow[r, "{\id \oplus r}"] & 
        \mathbb{C}^n \oplus \mathbb{C}^n \arrow[r, "\cong"] &
        \mathbb{C}^{2n} \arrow[r, "F_{2n}"] & {}\\
        \mathbb{C}^{2n} \arrow[r, "\cong"] & 
        \mathbb{C}^{n} \oplus \mathbb{C}^{n} \arrow[r, "\pi_1"] &
        \mathbb{C}^{n} \arrow[r, "\cdot \frac{1}{2}\omega_{2n}^{k / 2}"]& 
        \mathbb{C}^n &
    \end{tikzcd}
\end{center}

To ensure a cyclic distribution for the DFT, we take a cyclic distribution
of $\mathbb{C}^n$. Local index $k$ on $P(s)$ corresponds to global index
$s + kp$. If we embed this in the second component of
$\mathbb{C}^n \oplus \mathbb{C}^n \cong \mathbb{C}^{2n}$, we get global index
$n + s + kp = s + (n / p + k)p$, which shows it is the $k$th local index
in the second component of $\mathbb{C}^{n/p} \oplus \mathbb{C}^{n/p}$ on
$P(s)$. So the duplication $\langle \id, \id \rangle$ becomes a local
duplication under the cyclic distribution. As this is a local linear
computation, it fits in our framework as a computation superstep.

The function $\id \oplus r$ is not local under the cyclic distribution, so
we model it as a communication superstep. Local index $k$ on $P(s)$ has global
index $s + kp$, and 
$r(s + kp) = n - 1 - (s + kp) = p - 1 - s + (n / p - 1 - k)p$. So the reverse
is stored on $P(p - 1 - s)$ in local index $n / p - 1 - k$. To deal with the
embedding $\mathbb{C}^{n/p} \oplus \mathbb{C}^{n/p} \cong \mathbb{C}^{2n/p}$,
we first subtract $n/p$ from $k$ to obtain $k'  = k - n/p$, 
then reverse to give the index $n/p-k'-1$, and finally add $n/p$ back to the result, which gives local
index $3n / p - 1 - k$. In the notation of \cref{thm:tensor-of-bsp},
this gives the following communication step.

\[
    \rho(s, k) = \begin{cases}
                    k              & 0     \leq k <  n / p \\
                    3n / p - 1 - k & n / p \leq k < 2n / p \\
                 \end{cases}
\]
\[
    \psi(s, k) = \begin{cases}
                    s              & 0     \leq k <  n / p \\
                    p - 1 - s      & n / p \leq k < 2n / p \\
                 \end{cases}
\]

\subsubsection{Tensor Product}

To apply \cref{thm:tensor-of-bsp}, we must form the tensor product of
computation supersteps
\[
    \mathbb{C}^{n / p} \xrightarrow{\langle \id, \id \rangle}
    \mathbb{C}^{n / p} \oplus \mathbb{C}^{n / p} \cong \mathbb{C}^{2n/p}
\]
and
\[
    \mathbb{C}^{2n/p} \cong \mathbb{C}^{n / p} \oplus \mathbb{C}^{n / p} 
    \xrightarrow{\pi_1} \mathbb{C}^{n / p} 
    \xrightarrow{\cdot \frac{1}{2} \omega_{2n_l}^{k_l/2}} \mathbb{C}^{n/p}.
\]
The tensor product respects the product in the obvious way, e.g.
\[
    \bigotimes\limits_{l = 1}^d \langle \id_{c_l} \mid c_l \in [2] \rangle \cong
    \langle \id_{c_1} \otimes \cdots \otimes \id_{c_d} \mid c 
        \in [2] \times \cdots \times [2] \rangle.
\]
The projections are similarly combined dimension-wise.
To compose with \\
$\mathbb{C}^{2n/p} \cong \mathbb{C}^{n / p} \oplus \mathbb{C}^{n / p}$,
we can add $cn/p$ to the local index. The pointwise product is multiplied
across dimensions. This yields \cref{alg:dct:tensor}.

\begin{algorithm}[ht]
\caption{Parallel DCT-II for processor $P(s)=P(s_1, \cdots, s_d)$,
         rank-$d$}\label{alg:dct:tensor}
\begin{algorithmic}[1]
\Require $X: $ array of size $n_1 \times \cdots \times n_d$, distr$(X) = $ 
         rank-$d$ cyclic over $p_1 \times \cdots \times p_d$ processors such 
         that $p_l^2 | 2n_l$, for $l=1,\cdots,d$.
\Ensure $Y: $ array of size $n_1 \times \cdots \times n_d$, 
        distr$(Y) = $  rank-$d$ cyclic, such that 
        $Y = (C_{n_1} \otimes \cdots \otimes C_{n_d})(X)$.
    \For{$c \in [2] \times \cdots \times [2]$}
        \For{$k \in [n_1 / p_1] \times \cdots \times [n_d / p_d]$}
            \State $W^{(s)}[c n / p + k] = X^{(s)}[k]$;
        \EndFor
    \EndFor
    \For{$k \in [2n_1 / p_1] \times \cdots \times [2n_d / p_d]$}
        \State $W'^{((\psi_1 \times \cdots \times \psi_d)(s, k))}
            [(\rho_1 \times \cdots \times \rho_d)(k)] = W^{(s)}[k]$;
    \EndFor
    \State Compute $Z = DFT(W')$ with the parallel four-step FFT framework;
    \For{$k \in [n_1 / p_1] \times \cdots \times [n_d / p_d]$}
        \State $Y^{(s)}[k] = \prod_{l = 1}^d \omega_{2n_l}^{k_l / 2} Z^{(s)}[k]$;
	\EndFor
\end{algorithmic}
\end{algorithm}

\section{Related Work}

Johnson~\cite{Johnson90} has paved the way towards
a common framework for many DFT algorithms by using matrix decomposition.
This is also the fundamental unifying approach in the book by Van Loan~\cite{Loan92}.
The main difference with our work is that a matrix 
representation requires bases to be ordered linearly,
losing the higher-ranked structure of the DFT. This has led to a more
mechanical approach, deriving matrix-decompositions that are technically 
correct, but can be conceptually misleading.
For example, Johnson uses $F_{n / p} \otimes I_p$ for a cyclic distribution 
and $I_p \otimes F_{n / p}$ for a block distribution, while we would consider 
both operations as $\dsum F_{n / p}$ with a different isomorphism
$\dsum \mathbb{C}^{n / p} \cong \mathbb{C}^n$. The multiple $F_{n/p}$ in
\cref{chp1:four step sequential} arise from recursively applying the 
Discrete Fourier Transform on subarrays, 
not from multilinearity, which is why we believe the direct sum is a
better choice than the $\otimes$-notation.

Designing and implementing algorithms often goes hand-in-hand. We speculate
that previous work has sacrificed the higher-ranked structure because arrays
of arbitrary rank are challenging to work with in programming languages commonly
used for high-performance computing, such as Fortran, C,
 or C++. This is unfortunate,
as the structure of \cref{chp1:four step parallel:tensor}
does not depend on rank, so we would prefer to have one subroutine that
can take an array of any rank as argument. This is called 
\emph{rank-polymorphism} and is supported by languages like
APL~\cite{Iverson62} and Single-assignment~C~\cite{Scholz94}. We show that
rank-polymorphism can be used to describe tensor products of functions,
but it has also found applications in controlling concurrency and
data movement~\cite{Scholz23}. 

\section{Conclusion}

We introduced the notion of linear BSP algorithms: a class of parallel algorithms
that respects linearity and that follows  the Bulk Synchronous Parallel model for designing
algorithms. This model distributes data structures over multiple processors
with their own memory. Computation is done on these local data structures, and
separated from communication by global barrier synchronisation. The split of
computation and communication makes it possible to view algorithms as a
form of function decomposition into distributions, computation supersteps,
and communication supersteps. Investigating what the distributions,
computation, and communication look like in the category Vect led to the
main result of this paper: that linear BSP algorithms for a collection of
functions $f_1, \cdots, f_d$ give a linear BSP algorithm of the same
structure for their tensor product $f_1 \otimes \cdots \otimes f_d$
(\cref{thm:tensor-of-bsp}).

\section{Future Work}\label{subsec:future-tensor}

Future work could investigate what parallelism looks like in other categories,
and if it interacts with universal constructions in interesting ways.
For example, there are many algorithms that use complex conjugation,
which is not a linear function as $\overline{(\lambda \cdot z)} =
\overline{\lambda} \cdot \overline{z}$. The concept of conjugation
is captured and generalized in so-called $*$-algebras.
Functions that respect addition, but conjugate scalar multiplication are called
\emph{anti-linear}. These functions are ubiquitous in the study of Hilbert spaces,
making related categories of practical interest.

More efficient parallel algorithms for the DCT-II exist, based on using the DFT and 
functions that are linear with respect to $\mathbb{R}$, but not
$\mathbb{C}$~\cite{Narasimha78,Inda00}. This mix of ground fields keeps us from applying  
\cref{thm:tensor-of-bsp}, as we deal with $\otimes_{\mathbb{R}}$ and
$\otimes_{\mathbb{C}}$, which live in different categories. 
Future work could find a more efficient algorithm by parallelising Feig
and Winograd's algorithm~\cite{Feig92} using the BSP model. This algorithm
uses only real numbers, making it suitable for our approach.

Other interesting applications
that connect algebraic structures with applications can be found in
cryptography and topological data analysis.

The tensor product also describes the combination of Hadamard gates, making
it fundamental in quantum computing~\cite{Nielsen10}. As the performance of 
quantum computers relies on parallelism, it may also be a good candidate for
the higher level of abstraction that category theory can provide.

\bibliographystyle{abbrvurl}
\bibliography{paper.bib}

\end{document}